\documentclass{article}
\usepackage{etex}
\usepackage[usenames,dvipsnames,svgnames,table]{xcolor}
\usepackage{color}
\usepackage[latin1]{inputenc}
\usepackage{amsmath, enumerate, amsfonts}
\usepackage[widemargins]{a4}
\usepackage{tikz} 
\usepackage{tikz-cd}  
\usepackage{amsmath, enumerate, amsfonts,amssymb, amscd}

\newcommand{\pcom}{^\wedge_p}
\usepackage[all]{xy}
\xyoption{arc}
\newtheorem{Conjecture}{Conjecture}[section]

\newtheorem{Theorem}[Conjecture]{Theorem}

\newtheorem{Remark}[Conjecture]{Remark}

\usepackage[latin1]{inputenc}
\usepackage{graphicx}
\usepackage{amsmath, enumerate, amsfonts}
\usepackage[widemargins]{a4}
\usepackage[all]{xy}
\xyoption{arc}
\color{black}

\title{A few examples of $R-$good and $R-$bad classifying spaces}
\author{Nora Seeliger}
\date{}
\begin{document}
\maketitle
\begin{abstract}
For a commutative ring $R$, in contrast to the completion in the sense of Bousfield and Kan at just a prime integer, there cannot exist spaces which are good and bad in an arbitrary way.
\end{abstract}
\section{Introduction}
In general, it is a very hard question to decide whether a space is good or bad in the sense of Bousfield and Kan [\textit{Homotopy Limits, Completions and Localizations}, Springer Lecture Notes in Mathematics, Springer-Verlag, Berlin, Heidelberg, New York, 1972] and there are not many explicit examples in the literature. In this note we give results about $R-$completed classifying spaces for various solid rings $R$. In particular, we show that spaces cannot be good and bad at different solid rings in any arbitrary combination. The special case for primes only was a question raised during the talk of Robert Oliver on $p-$local homotopy theory at the MSRI Berkeley at the beginning of the algebraic topology program in 2014 concerning the possible compatibilty of completions at two different primes at the same time and was answered positively by the author in \cite{pgoodpbad} 
who was supported by a postdoctoral fellowship at the University of Haifa under the direction of Professor David Blanc. I would like to thank him and Emmanuel Dror Farjoun for discussions on the topic.
\section{Preliminaries}
The $p$-completion $X\pcom$ of a topological space $X$ in the sense of Bousfield-Kan \cite{BK} is a special case of the $R-$completion functor of Bousfield and Kan which is defined for solid rings $R$. This is a functor from the category of simplicial sets to itself. The category of simplicial sets will be denoted by spaces thoughout this article. A ring $R$ is \textbf{solid} if the map $R\otimes _{\mathbb{Z}}R\rightarrow R$ is an isomorphism. The solid rings are $\mathbb{Z}/n\mathbb{Z}$, subrings of the rationals $\mathbb{Z}[J^{-1}]$ for any set $J$ of primes, the product rings $\mathbb{Z}[J^{-1}]\times \mathbb{Z}/n\mathbb{Z}$ where each prime factor of $n$ is in $J$, and direct limits of these three types of rings \cite{BK2}. The Bousfield-Kan completion is related to completions and localizations in the sense of Malcev \cite{Malcev}, Sullivan \cite{Sullivan}, and Quillen \cite{Quillen}.
 In $p-$local homotopy theory the
classifying space  $|\mathcal{L}|\pcom$ of a $p-$local finite group $(S,\mathcal{F},\mathcal{L})$ is one of the main 
objects of study where $(-)\pcom$
denotes the $p-$completion functor with respect to the ring $\mathbb{Z}/p\mathbb{Z}$.  A space X is called \textbf{$R-$complete} if the map $X\rightarrow X^{\wedge}_R$ is a mod $R-$equivalence and  \textbf{$R-$good} if the natural map $(H_*;R)\rightarrow H_*(X^{\wedge}_R ;R)$ is an isomorphism. Otherwise the space is called \textbf{$R$-bad}. Examples of spaces which are good for all primes and the integers are classifying spaces of finite groups. Spaces with finite homotopy groups are good for all solid rings $R$. A space $X$ is called $p-$good for a prime $p$ if it is good for the ring $\mathbb{Z}/p\mathbb{Z}$. For spaces with finite homotopy groups or finite homology groups the $\mathbb{Z}-$completion is up to homotopy, the product of the $\mathbb{Z}/p\mathbb{Z}-$completions. The space $\mathbb{R}P^2$ is good for all primes and bad for $\mathbb{Z}_{(J)}$ as long as $2\in J$, and bad for the integers \cite{BK}. Recall that the space $S^1\vee S^1$ is bad for all primes and the integers \cite{Bousfield}.
 A group $\pi$ acts on a group $G$
if there is given a homomorphism
$\alpha :\pi\rightarrow Aut(G)$. Such an action is called \textbf{nilpotent} if there exists a finite
sequence of subgroups of $G$:
$G=G_1\supset\cdots\supset G_j\supset\cdots\supset G_n=1$
such that for each $j$ we have
$G_j$ is closed under the action of $\pi$,
$G_{j +1}$ is normal in $G_j$, and $G_j/G_{j +1}$ is abelian, and
the induced action on $ G_j/G_{j+1} $ is trivial.\\
A (possibly infinite) group $G$ has a \textbf{Sylow} $p-$subgroup $S$ if it has a $p-$subgroup isomorphic to $S$ and all other subgroups of $G$ whose order is a power of $p$ are subconjugate to $S$. Examples of infinite groups with finite Sylow $p-$subgroups are group models for fusion systems in the sense of Robinson and Leary-Stancu and ourselves together with Leip and ourselves together with Libman. Subgroups of amalgamated products are described by Kurosh's subgroup Theorem. The most general version states that \cite{Massey} if $H$ is a subgroup of the free amalgamated product of groups $ \underset{i\in I}{*}G_i = G$, then $H=F(X)*(*_{j\in J} g_jH_jg_j^{-1})$, where $X$ is a subset  of $G$ and $J$ is some index set and $g_j \in G$ and each $H_j$ is a subgroup of some $G_i$. Recall that for any pair of groups $ G,H$  we have a weak equivalence of classifying spaces $B(G*H)\simeq BG\vee BH$, \cite{Massey}. The mod $R-$Fibre Lemma \cite[Lemma 5.1]{BK} states that the $R$-completion preserves fibrations of connected spaces up to homotopy $F\rightarrow E\rightarrow B$ for which the fundamental group of the base space $\pi _1(B)$ acts nilpotently on every reduced homology group of the fibre with coefficients in the ring of the completion $\overline{H}_i(F;R)$. It is used to show that spaces are good or bad without constructing the completion tower. A group $G$ is \textbf{virtually finite} if it has a free subgroup $H<G$ of finite index. A group $G$ is \textbf{a virtual finite $p$-group} if it has a free subgroup of index a power of $p$. A space is \textbf{nilpotent} if the fundamental group acts nilpotently on all higher homotopy groups \cite{BK}. For a ring $R$ a group $G$ is $R-$nilpotent if it has a finite central series $G=G_1\subset\cdots G_j\cdots G_n=1 $ such that each quotient $G_j/G_{j+1}$ admit an $R$-module structure.A map $f:\{G_s\}\rightarrow \{H_s\}$ between two towers of groups is a \textbf{pro-isomorphism} if, for every group $B$, it induces an isomorphism 
$\underset{\leftarrow}{lim}\text{ }Hom_{(groups)}(H_s,B)\cong \underset{\leftarrow }{lim}\text{ }Hom_{(groups)}(G_s,B)$.
The completion tower $R_nX$defined by Bousfield and Kan preserves the homology with $R$-coefficients \cite{Farjoun2}. This property of preserving $R$-homology characterizes the tower completely \cite{BK}. 
For a ring $R$ and a space $X$ denote $E^RX$  the $HR-$localization functor defined in \cite{Bousfield1}.
While it is possible to prove that a space is bad for prime numbers using the Fibre Lemma, this is not possible for arbitrary solid rings as we showed in \cite{pgoodpbad}.\\[0.3cm] The following proves that in contrast to the case of primes in the general case a space cannot be both $R$-good and $R'$-bad for arbitrary solid rings $R$ and $R'$ at the same time. 
\begin{Theorem}
Let either
\begin{itemize}
\item $R=\mathbb{Z}/m\mathbb{Z}$ and $R'=\mathbb{Z}/n\mathbb{Z}$, where $m$ and $n$ are integers such that $m$ and $n$ are coprime, or
\item $R=\mathbb{Z}[J^{-1}]$ and $R'=\mathbb{Z}/n\mathbb{Z}$, where $n$ is an integer all of whose  prime factors are in $J$.
\end{itemize}

Let $X$ be a space which is good for the ring $R$ and bad for the ring $R'$. Then $X$ is bad for $R\times R'$.
\end{Theorem}
\underline{Proof:} We have from \cite[Proposition 9.5]{BK} that $X\widehat{_{R\times R'}}\simeq X^{\wedge}_R\times X^{\wedge}_{R'}$. The result then follows from the K\"unneth Theorem.$\Box$
\begin{Remark}
Let $p$ be a prime and $X$ a $p-$good space. Then $X$ can be both good or bad for $\mathbb{Z}_{(p)}$.
\end{Remark}
\begin{Theorem}[\cite{IvanovMikhailov1}{, Corollary 9.2.}]
For a finitely presented metabelian group $X$, the natural map $E^RX\rightarrow \hat{X}_R$ is an isomorphism for $R=\mathbb{Q}$ or $\mathbb{Z}/n\mathbb{Z}$.
\end{Theorem}
\begin{Remark}[\cite{IvanovMikhailov1}{, Theorem 9.1.}]
The above result cannot be generalized to the case $R=\mathbb{Z}$.
\end{Remark}
\begin{Remark}
Recall that $\mathbb{R}P^{\infty}\simeq BC_2$ is good for the integers and all primes but bad for $\mathbb{Z}_{(2)}$. 
\end{Remark}
\begin{Remark}
The space $\mathbb{R}P^2$ is bad for the prime $2$ and good for all other primes because the reduced homology $\overline{H}(\mathbb{R}P^2;\mathbb{F}_p)$ is trivial for $p$ odd.
\end{Remark}

\begin{Theorem}
Let $R$ be a nonzero solid ring. The classifying space $S^1\vee S^1\simeq B(\mathbb{Z}*\mathbb{Z})$ is $R$-bad.
\end{Theorem}
\underline{Proof:} Assume the contrary that $S^1\vee S^1$ is $R-$good. Then any finite wedge is $R-$good as well. Recall that $(S^1\vee S^1)^{\wedge}_R\simeq (B(\mathbb{Z}*\mathbb{Z}))^{\wedge}_R\simeq B((\mathbb{Z}*\mathbb{Z})^{\wedge}_R)$ in this situation because $\mathbb{Z}*\mathbb{Z}$ has finite virtual cohomological dimension. This is a contradiction in low dimensions because through the group completion we have that $(\mathbb{Z}*\mathbb{Z})^{\wedge}_R$ is not finitely generated in low degrees as proved in \cite{Bousfield1}.$\Box$

Dr.\ Nora Seeliger, Department of Mathematics, University of Haifa, Faculty of Natural Sciences, Science and Education Building, Room 615, Abba Khoushy Avenue 199, 3498838 Haifa, ISRAEL.\\Email:
nseelige@campus.haifa.ac.il

\begin{thebibliography}
{}
\bibitem{AKO} M.\ Aschbacher, R.\ Kessar, B.\ Oliver, \textit{Fusion systems in algebra and topology}, London Mathematical Society Lecture Note Series: 31, Cambridge University Press, 2011.
\bibitem{Bousfield1} A.\ K.\ Bousfield, \textit{Homological localization towers for groups and $\pi -$modules}, Mem.\ Amer.\ Math.\ Soc, no. 186, 1977. 
\bibitem{Bousfield} A.\ K.\ Bousfield, \textit{On the $p-$adic completions of nonnilpotent spaces}, Transactions of the American Mathematical Society, Volume 331, Number 1, May 1992, 335--359.
\bibitem{BK} A.\ K.\ Bousfield, D.\ M.\ Kan, \textit{Homotopy Limits, Completions and Localizations}, Springer Lecture Notes in Mathematics, Springer-Verlag, Berlin, Heidelberg, New York, 1972.
\bibitem{BK2} A.\ K.\ Bousfield, D.\ M.\ Kan, \textit{The core of a ring}, J.\ Pure Applied Algebra 2 (1972), 73--81.
\bibitem{CurtisReiner} C.\ W.\ Curtis, I.\ Reiner, \textit{Representation Theory of Finite Groups and Associative Algebras}, Interscience Publishers, New York, 1962.
\bibitem{DwyerFarjounRavenel} W.\ G.\ Dwyer, E.\ D.\ Farjoun, and D.\ C. Ravenel, \textit{Bousfield localizations of classifying spaces of nilpotent groups}, Proceedings of the American Mathematical Society, Volume 127, Pages 1855--1861.
\bibitem{Farjoun} E.\ D.\ Farjoun, \textit{Two Completion Towers for Generalized Homology}, Contemporary Mathematics Volume 265 (2000).
\bibitem{Farjoun2} E.\ D.\ Farjoun, \textit{Pro-nilpotent representation of homology types}, Proceedings of the American Mathematical Society, Volume 38, Number 3, May 1973, 657--660.
\bibitem{Hall} M. Hall, \textit{The theory of groups}, Macmillan (1959).

\bibitem{IvanovMikhailov1} S.\ O.\ Ivanov, R.\ Mikhailov, \textit{On a problem of Bousfield for metabelian groups}, Adv.\ Math.\ 290 (2016), 552--589. 
\bibitem{IvanovMikhailov2} S.\ O.\ Ivanov, R.\ Mikhailov, \textit{On lenghts of $H\mathbb{Z}-$localization towers}, preprint.





\bibitem{Malcev} A.\ L.\ Malcev, \textit{Nilpotent groups without torsion}, Jzv. Akad. Nauk. SSSR
, Math. 13 (1949), 201--212.
\bibitem{MartinoPriddy} J.\ Martino, S.\ Priddy, \textit{Unstable homotopy classification of $BG\pcom$}, Math.\ Proc.\ Camb.\ Phil.\ Soc.\ 137 (2004) 321--347.
\bibitem{Massey} William S. Massey, \textit{Algebraic topology: an introduction}, Graduate Texts in Mathematics, Springer-Verlag, New York, 1977.
%
\bibitem{Quillen} D.\ G.\
 Quillen, \textit{Rational homotopy theory
 }, Annals of Math. 90
 (1969), 205--295.
\bibitem{pgoodpbad} N.\ Seeliger, \textit{A few examples of $p$-good and $p$-bad classifying spaces}, 	arXiv:1411.7490.
\bibitem{Sullivan} D.\ Sullivan, \textit{Geometric topology, part I: localization, periodicity and Galois symmetry}, MIT (1970).
\end{thebibliography}
\end{document}